# FLUID LIMITS OF PURE JUMP MARKOV PROCESSES: A PRACTICAL GUIDE

**R. W. R. Darling**

National Security Agency
P. O. Box 535, Annapolis Junction, MD 20701
`rwrd@afterlife.ncsc.mil`

July 10, 2002

ABSTRACT: A rescaled Markov chain converges uniformly in probability to the solution of an ordinary differential equation, under carefully specified assumptions. The presentation is much simpler than those in the outside literature. The result may be used to build parsimonious models of large random or pseudo–random systems.

## 1. INTRODUCTION

■ **1.1 Our Goal**

In many fields of research one seeks to give a parsimonious description of the behavior of a large system whose dynamics depend on the random interactions of many components (molecules, genes, consumers, …). Such a description may take of the form of the solution of an ordinary differential equation, derived as the deterministic limit of a suitably scaled Markov process, as some scale parameter $N \longrightarrow \infty$. We are especially interested in applications to large random combinatorial structures, in order to extend results such as those presented in Bollobás (2001), Kolchin (1999), and Janson, Łuczak and Rucinski (2000).

The study of such limit theorems is a rich, and technically advanced, topic in probability: see the books of Ethier and Kurtz (1986) (esp. p. 456), and Jacod and Shiryaev (1987) (esp. p. 517). The work of Aldous, e.g. Aldous (1997), offers many examples of the application of modern probabilistic limit arguments to random combinatorial structures.

Our goal here is to establish a fairly general fluid limit result, accessible to anyone familiar with basic Markov chains, Poisson processes, and martingales. For a more powerful result using Laplace transforms of Lévy kernels, see Darling and Norris (2001). One example is worked out in detail; others can be found in Darling and Norris (2001).

### ■ 1.2 Prototype of the Theorem We Seek

Readers will have encountered the Weak Law of Large Numbers (WLLN) for the sum

$$X_n \equiv \sum_{i=1}^{n} U_i$$

of independent, identically distributed random variables $(U_i)_{i \geq 1}$ with common mean $\mu$ and variance $\sigma^2$; namely that $N^{-1} X_N$ converges in probability to $\mu$. A more sophisticated technique is to apply a maximal inequality for submartingales (Kallenberg (2002), p. 128) to the square of the martingale

$$(N^{-1}(X_n - n\mu))_{0 \leq n \leq N} \tag{1}$$

to obtain:

$$\delta^2 \mathbb{P}\left[\max_{n \leq N} \left(\frac{X_n - n\mu}{N}\right)^2 \geq \delta^2\right] \leq \mathbb{E}\left[\left(\frac{X_N - N\mu}{N}\right)^2\right] = \frac{\sigma^2}{N}.$$

This is a functional form of the WLLN for the rescaled Markov chain $Y_t^N \equiv N^{-1} X_{Nt}$, for $t$ a multiple of $\epsilon \equiv N^{-1}$, namely

$$\mathbb{P}\left[\max_{t \leq 1} |Y_t^N - t\mu| \geq \delta\right] = O(N^{-1}).$$

Our aim is to generalize such a fluid limit result from random walks, as above, to a suitable class of Markov chains in $E$, where $E$ is a Euclidean space or a separable Hilbert space. The fluid limit will in general not be of the form $y[t] = t\mu$, but the unique solution of some ordinary differential equation $\dot{y}[t] = b[y[t]]$ in a suitably selected domain $S \subset E$. The first exit time from $S$ will almost always be an important consideration. Indeed the most common application of such a theorem is to prove:

*The first exit time from S of the rescaled Markov chain converges in probability to the first exit time of the fluid limit.*

## ■ 1.3 Illustration

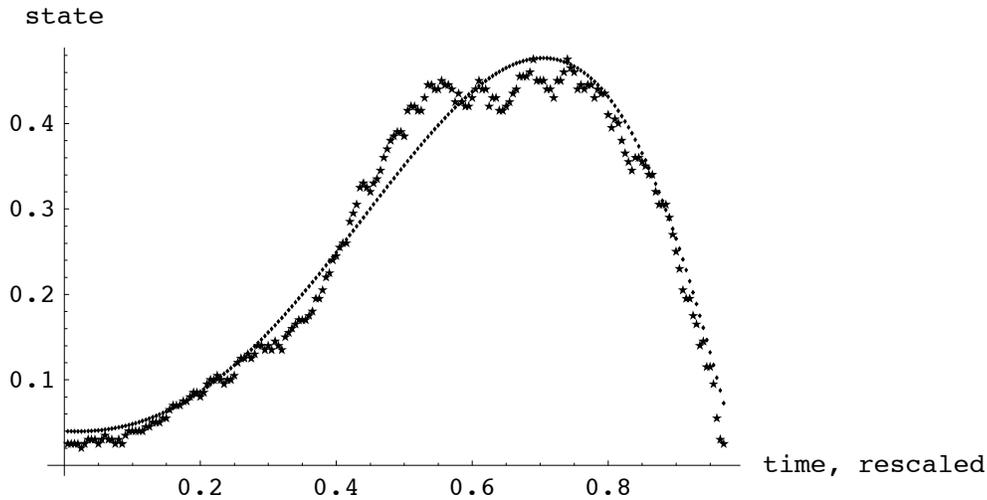

Scaling parameter $N = 200$: $\star$ tracks a rescaled Markov chain — the patch process in Darling and Norris (2001) — and $\diamond$ tracks its fluid limit, whose derivation is explained below. We are interested in bounds, in terms of $N$, on the maximum distance between the two trajectories.

# 2. PURE JUMP MARKOV PROCESSES

## ■ 2.1 Elementary Construction of Pure Jump Processes

The reader is assumed to be familiar with discrete–time Markov chains and with Poisson processes. We shall give here a naïve presentation of Markov processes of pure jump type, sufficient for our needs.

Suppose $(X_n)_{n \geq 0}$ is a discrete–time Markov chain, whose state space is some subset $I \subseteq E$. We assume that $X_0$ and the increments of the chain have finite means and covariances:

$$\mathbb{E}[X_{n+1} - X_n \mid X_n = x] = \mu[x]; \tag{2}$$

$$\mathrm{Var}[X_{n+1} - X_n \mid X_n = x] = \Sigma[x]. \tag{3}$$

Let $(v[t])_{t \geq 0}$ be a Poisson process, with event times $\tau_1 < \tau_2 < \ldots$, which is dependent on $(X_n)_{n \geq 0}$ in the following sense: there is some bounded **rate function** $c : I \longrightarrow (0, \infty)$, such that each inter–event time $\tau_{n+1} - \tau_n$ is Exponential with mean $1/c[x]$ on the event $\{X_n = x\}$. This implies that

$$\text{Var}[\tau_{n+1} - \tau_n \mid X_n = x] = c[x]^{-2}. \tag{4}$$

As in Kallenberg (2002), Chapter 12, the formula

$$Y_t \equiv X_{\nu[t]} \tag{5}$$

defines a **pure jump Markov process** $(Y_t)_{t \geq 0}$. In other words $Y_{\tau_n} = X_n$, so $(Y_t)_{t \geq 0}$ has the same increments as does $(X_n)_{n \geq 0}$, and these occur at the random times $\tau_1 < \tau_2 < \ldots$. The effect is to transform a discrete–time Markov chain into a continuous–time Markov process with a possibly variable jump rate.

### ■ 2.2 Why the Random Jump Times?

The reader may be wondering at this point why we focus hereafter on the continuous–time process $(Y_t)_{t \geq 0}$, rather than on the apparently simpler Markov chain $(X_n)_{n \geq 0}$. There are two reasons:

• We have to make the transition from discrete time to continuous time at some place in the argument, and probabilists find it convenient to do this at the very beginning.

• There are many models in which we want the jump rate to vary according to some function of the process; the construction above allows this.

### ■ 2.3 An Associated Martingale

Define $b[x] \equiv c[x]\mu[x]$, and define

$$A_t \equiv Y_0 + \int_0^t b[Y_{s_-}]\,ds.$$

The technical term for $(A_t)_{t \geq 0}$ is the **compensator** of $(Y_t)_{t \geq 0}$; see Kallenberg (2002), Chapter 25. A central technical role, comparable to that of the martingale (1), is played by the difference between the process and its compensator, namely

$$M_t \equiv Y_t - A_t. \tag{6}$$

By construction, $(M_t)_{t \geq 0}$ is a right continuous vector process with left limits, and with $M_0 = 0$, and its values at the jump times are given by:

$$M_{\tau_n} = X_n - X_0 - \sum_{j=1}^{n} (\tau_j - \tau_{j-1}) b[X_{j-1}], \tag{7}$$

*Fluid Limits* 5taking $\tau_0 \equiv 0$. In the modern approach to Markov processes (Kallenberg (2002), p. 382), the definition of the process $(Y_t)_{t \geq 0}$ renders $(M_t)_{t \geq 0}$ a martingale, under mild regularity conditions. Since we are not assuming this technical apparatus, we shall prove this directly.

### ■ 2.4 Proposition

*Assume that there exist constants $C_1$, $C_2$, such that*

$$\sup_{x \in I} \{\text{Trace}[\Sigma[x]] + \|\mu[x]\|^2\} \leq C_1; \tag{8}$$

$$\sup_{x \in I} c[x] \leq C_2. \tag{9}$$

*Then $(M_t)_{t \geq 0}$ is a square integrable vector martingale, and for each $u > 0$, and each $\delta > 0$,*

$$\mathbb{P}[\sup_{t \leq u} \|M_t\|^2 \geq \delta] \leq \min_{n > 1} \left\{ \mathbb{P}[\xi_n < u] + \frac{nC_1}{\delta} \right\}. \tag{10}$$

*where $\xi_n$ is a Gamma$(n, 1/C_2)$ random variable.*

The proof, which is not difficult, is postponed to Section 3.1.

### ■ 2.5 A Sequence of Pure Jump Processes

Suppose that for each $N \geq 1$ we have a Markov chain $(X_n^N)_{n \geq 0}$, on a state space $I_N \subseteq E$, whose increments have mean and covariance given by $\mu_N[x]$ and $\Sigma_N[x]$, respectively, and a rate function $c_N[x]$. As above, we define $b_N[x] \equiv c_N[x]\mu_N[x]$. Following the construction (5), we obtain a pure jump Markov process $(Y_t^N)_{t \geq 0}$ for each $N$, with associated compensator $(A_t^N)_{t \geq 0}$ and martingale $(M_t^N)_{t \geq 0}$.

Let $D \subseteq E$ be any closed set such that $D \supseteq \bigcup I_N$. We fix a relatively open set $S \subseteq D$, and define $S_N \equiv S \cap I_N$. The set $S$ is the region in space where laws of the processes will converge, but we make slight use of $D$ in the proof of Theorem 2.8.

We assume that the parameters of these processes are convergent in the following way: here $\kappa_1[\delta]$, $\kappa_2$, $\kappa_3$ denote positive constants, and the inequalities hold uniformly in $N$.

- **Initial conditions** converge: assume that there is some $a \in \bar{S}$, the relative closure of $S$ in $D$, such that, for each $\delta > 0$,

$$\mathbb{P}[\|Y_0^N - a\| > \delta] \leq \kappa_1[\delta]/N. \tag{11}$$

- **Mean dynamics** converge: assume that there is a Lipschitz vector field $b : D \longrightarrow E$ (the values of $b$ outside $\bar{S}$ are irrelevant) such that

*R. W. R. Darling* Printed by Mathematica for Students *July 2002*

$$\sup_{x \in S_N} \|b_N[x] - b[x]\| \longrightarrow 0. \tag{12}$$

- **Noise converges to zero**:

$$\sup_{x \in S_N} c_N[x] \leq \kappa_2 N; \tag{13}$$

$$\sup_{x \in S_N} \{\text{Trace}[\Sigma_N[x]] + \|\mu_N[x]\|^2\} \leq \kappa_3 N^{-2}. \tag{14}$$

In effect we are choosing "hydrodynamic scaling": the increments of $(X_n^N)_{n \geq 0}$ are $O(N^{-1})$, and the jump rate is $O(N)$. This is not the only possibility, but is consistent with (12). The purpose of (14) and (13) is to control the martingale part of the processes, as we see in (26).

### ■ 2.6 How to Specify *D* and *S* in Practice

In practical cases, one will usually identify first the formula for the limiting vector field $b[x]$. One will then identify a set on which it is Lipschitz (for example, by studying where its partial derivatives are uniformly bounded), and choose $S$ within this set. Then take $D$ to be a convenient superset.

### ■ 2.7 Fluid Limit

Since $b$ is Lipschitz on $\overline{S} \subseteq D$, there is a unique solution $(y[t])_{0 \leq t \leq \zeta[a]}$ in $\overline{S}$, where $\zeta[a] \equiv \inf \{t \geq 0 : y[t] \notin S\} \leq \infty$, to the ordinary differential equation:

$$\dot{y}[t] = b[y[t]], \quad y[0] = a. \tag{15}$$

We propose to show that this solution is the **fluid limit** of the sequence of Markov processes, in the following sense:

### ■ 2.8 Theorem

*Assume (11), (12), (13) and (14). Let $\sigma_N \equiv \inf \{t \geq 0 : Y_t^N \notin S\}$. Fix $\delta > 0$.*

*(i) For any finite time $u \leq \zeta[a]$,*

$$\mathbb{P}[\sup_{0 \leq t \leq u} \|Y_{t \wedge \sigma_N}^N - y[t \wedge \sigma_N]\| > \delta] = O(N^{-1}). \tag{16}$$

*(ii) Suppose now that $a \in S$. Then for any finite time $u < \zeta[a]$, $\mathbb{P}[\sigma_N < u] = O(N^{-1})$; hence for such $u$,*

$$\mathbb{P}[\sup_{0 \leq t \leq u} \|Y_t^N - y[t]\| > \delta] = O(N^{-1}). \tag{17}$$

*Moreover if $\zeta[a] < \infty$ (i.e. the ODE solution leaves S in finite time), then*

$$\mathbb{P}[\sup_{0 \le t \le \zeta[a]} \|Y^N_{t \wedge \sigma_N} - y[t]\| > \delta] = O(N^{-1}); \quad (18)$$

$$\mathbb{P}[|\sigma_N - \zeta[a]| > \delta] = O(N^{-1}). \quad (19)$$

### 2.8.1 Remarks:

• Frequently one encounters examples in which the obvious vector field $b$ is not globally Lipschitz. However for a finite $u > 0$, there may be a unique solution to $\dot{y}[t] = b[y[t]]$ for $t \in [0, u]$ started at $y[0] = a$, and the restriction of $b$ to a closed set $D \supseteq (y[t])_{0 \le t \le u}$ is Lipschitz. Apply the Theorem to any relatively open set $S \subseteq D$.

• The modification "$t \wedge \sigma_N$" in (18) is essential. For any $N$, there is a good chance that the Markov process will leave $S$ before time $\zeta[a]$, whereupon it may become subject to entirely different dynamics.

• The Markov property is used only to establish that $(M_t)_{t \ge 0}$ is a vector martingale. The proof works equally well if $(X_n)_{n \ge 0}$ is merely the **image** under a suitable mapping $\Phi$ of a Markov chain on an arbitrary measurable space, and the definition of compensator is adjusted accordingly. In that case $(X_n)_{n \ge 0}$ itself need not be Markov. We give an example below where $\Phi$ is linear.

## 3. PROOFS

### 3.1 Proof of Proposition 2.4

### 3.1.1 Step I

*Each jump time $\tau_n$ is bounded below by a Gamma$(n, 1/C_2)$ random variable $\xi_n$, in the sense that*

$$\mathbb{P}[\tau_n \le t] \le \mathbb{P}[\xi_n \le t] < (tC_2)^{n-2}/(n-1)!, \quad \forall\, t > 0, \quad \forall\, n > 1. \quad (20)$$

A consequence of (9) is that the $n$–th jump time, $\tau_n$, is minorized by $\xi_n$. Hence

$$\mathbb{P}[\tau_n \le t] \le \mathbb{P}[\xi_n \le t].$$

An elementary estimate, based on an inequality of the form

$$\int_0^t z^{n-1} e^{-z}\, dz < t^{n-2} \int_0^t z e^{-z}\, dz < t^{n-2},$$

gives (20).

■ **3.1.2 Step II**

*Assume $\mathbb{E}[\|X_0\|^2] < \infty$. The processes $(Y_t)_{t \geq 0}$, $(A_t)_{t \geq 0}$, $(M_t)_{t \geq 0}$ and $(M_{\tau_n})_{n \geq 0}$ are square integrable.*

From (8) we have $\mathbb{E}[\|X_{n+1} - X_n\|^2] \leq C_1$, so by Pythagoras and induction,

$$\mathbb{E}[\|X_n - X_0\|^2] \leq nC_1, \ n = 1, 2, \ldots. \tag{21}$$

Conditioning on the trajectory of the Markov chain $(X_n)_{n \geq 0}$,

$$\mathbb{E}[\|Y_t - Y_0\|^2 \mid (X_n)_{n \geq 0}] = \mathbb{E}\left[\sum_n 1_{\{\tau_n \leq t \leq \tau_{n+1}\}} \|X_n - X_0\|^2 \ \bigg| \ (X_n)_{n \geq 0}\right]$$

$$\leq \sum_n \|X_n - X_0\|^2 \mathbb{P}[\tau_n \leq t \mid (X_n)_{n \geq 0}] \leq \sum_n \|X_n - X_0\|^2 \mathbb{P}[\xi_n \leq t].$$

The expectation of this is finite by (20) and (21). Hence $(Y_t)_{t \geq 0}$ is square integrable.

It is immediate from (7) that

$$M_{\tau_{n+1}} - M_{\tau_n} = X_{n+1} - X_n - (\tau_{n+1} - \tau_n) b[X_n]. \tag{22}$$

From (2), (3), and (4), we see that the conditional distribution of (22), given $X_n = x$, has mean $\mu[x] - b[x]/c[x] = 0$, and covariance

$$\Sigma[x] + c[x]^{-2} b[x] b[x]^T = \Sigma[x] + \mu[x] \mu[x]^T.$$

From (8), it follows that

$$\sup_n \mathbb{E}[\|M_{\tau_{n+1}} - M_{\tau_n}\|^2] \leq C_1.$$

Apply Pythagoras and induction, and the fact that $M_0 = 0$, to deduce that, for all $n$,

$$\mathbb{E}[\|M_{\tau_n}\|^2] \leq nC_1. \tag{23}$$

A similar argument shows that $\mathbb{E}[\|A_{\tau_n} - A_0\|^2] \leq 2nC_1$. For $\{\tau_n \leq t < \tau_{n+1}\}$, $A_t$ interpolates linearly between $A_{\tau_n}$ and $A_{\tau_{n+1}}$. By an argument similar to that for $Y_t$,

$$\mathbb{E}[\|A_t - A_0\|^2] < \infty.$$

Finally the square integrability of $M_t$ follows because $M_t = Y_t - A_t$.

### 3.1.3 Step III

$(M_t, \mathcal{F}_t)_{t \geq 0}$ *is a square integrable vector martingale, where* $\mathcal{F}_t \equiv \sigma\{Y_s, 0 \leq s \leq t\}$.

We saw in Step II that each $M_t$ is square integrable. To prove the martingale property, i.e. that

$$\mathbb{E}[M_{u+t} \mid \mathcal{F}_u] = M_u, \ \forall \, u \geq 0, \ \forall \, t \geq 0,$$

it suffices by the Markov property to prove it when $u = 0$, i.e. that

$$\mathbb{E}[M_t] = 0, \ t \geq 0. \tag{24}$$

There is no loss of generality in assuming that $X_0 \equiv Y_0$ is non–random, and equals some specific $x$. Because $M_t = -tb[x]$ for all $t < \tau_1$, we may write

$$M_t = (M_{\tau_1} + (M_t - M_{\tau_1}))1_{\{t \geq \tau_1\}} - tb[x]1_{\{t < \tau_1\}}. \tag{25}$$

Abbreviate $c[x]$ to $\lambda$, so $\lambda\mu[x] = b[x]$. Since $\tau_1$ has density $f[z] \equiv 1_{\{z>0\}}\lambda e^{-\lambda z}$,

$$\mathbb{E}[\tau_1 1_{\{\tau_1 \leq t\}}] = \int_0^t z f[z] dz = \lambda^{-1}(1 - e^{-\lambda t}) - t e^{-\lambda t}.$$

From (7), taking $n = 1$,

$$\mathbb{E}[M_{\tau_1} 1_{\{t \geq \tau_1\}} - tb[x]1_{\{t < \tau_1\}}] = \mathbb{E}[(X_1 - x - \tau_1 b[x])1_{\{t \geq \tau_1\}}] - tb[x]e^{-\lambda t},$$

and an easy calculation based on the previous integral show that this is zero. From (25), bearing in mind that $M_t$ is integrable and that $\tau_0 = 0$, we have now proved that

$$\mathbb{E}[M_t] \equiv \mathbb{E}[(M_t - M_{\tau_0})1_{\{t \geq \tau_0\}}] = \mathbb{E}[(M_t - M_{\tau_1})1_{\{t \geq \tau_1\}}].$$

We can repeat the same steps $n$ times, to obtain

$$\mathbb{E}[M_t] = \mathbb{E}[(M_t - M_{\tau_n})1_{\{t \geq \tau_n\}}].$$

Apply Lebesgue Monotone Convergence and (20) to obtain:

$$\lim_{n \to \infty} \mathbb{E}[\|M_t\| 1_{\{t \geq \tau_n\}}] = 0.$$

Apply the Cauchy–Schwartz inequality, (20), and (23), to obtain

$$\mathbb{E}[\|M_{\tau_n}\| 1_{\{t \geq \tau_n\}}]^2 \leq \mathbb{E}[\|M_{\tau_n}\|^2] \mathbb{P}[\tau_n < t] \longrightarrow 0.$$

This establishes (24), as desired.

### 3.1.4 Step IV

Given $\delta > 0$, $u > 0$, and $n \geq 1$, we may write

$$\mathbb{P}[\sup_{t \leq u} \|M_t\|^2 \geq \delta] \leq \mathbb{P}[\tau_n < u] + \mathbb{P}[\sup_{t \leq \tau_n} \|M_t\|^2 \geq \delta].$$

Apply a standard maximal inequality (Kallenberg (2002), p. 128) to the submartingale $(\|M_{t \wedge t_n}\|)_{t \in \mathbb{Q}_+}$, to obtain, for any $u > 0$ and $\epsilon > 0$,

$$\mathbb{P}[\sup_{t \leq \tau_n} \|M_t\|^2 \geq \delta] \leq \delta^{-1} \mathbb{E}[\|M_{\tau_n}\|^2].$$

By (20) and (23),

$$\mathbb{P}[\sup_{t \leq u} \|M_t\|^2 \geq \delta] \leq \mathbb{P}[\xi_n < u] + \delta^{-1} n C_1.$$

This holds for every $n$, so we obtain (10). □

## 3.2 Proof of Theorem 2.8:

### 3.2.1 Step I

*Suppose (14) and (13) hold when $S \equiv D$. For each $u > 0$, and each $\delta > 0$, the martingale part of $(Y_t^N)_{t \geq 0}$ satisfies*

$$\mathbb{P}[\sup_{t \leq u} \|M_t^N\| \geq \delta] = O(N^{-1}). \tag{26}$$

**Proof:** The effect of taking $S \equiv D$ is that we do not have to worry about exit from $S$. We apply (10), taking

$$C_1 \equiv \kappa_3 N^{-2}; \quad C_2 \equiv \kappa_2 N; \quad n \equiv u \kappa N$$

for some $\kappa > \kappa_2$. This immediately implies

$$\delta^{-1} n C_1 = O(N^{-1}).$$

The Gamma random variable $\xi_n$ in (10) has mean strictly greater than $u$, and variance which is $O(N^{-1})$. By Chebyshev's Inequality,

$$\mathbb{P}[\xi_n < u] = O(N^{-1}).$$

The assertion (26) follows. □

### 3.2.2 Step II

First we consider the case where $S \equiv D$. On taking the difference between the two equations

$$y[t] = a + \int_0^t b[y[s]]ds; \tag{27}$$

$$Y_t^N = Y_0^N + \int_0^t b_N[Y_s^N]ds + M_t^N$$

(we may replace $Y_{s-}^N$ by $Y_s^N$ in the integrand), we find that

$$\|Y_t^N - y[t]\| \leq \|Y_0^N - a\| + \|M_t^N\| + \int_0^t \|b_N[Y_s^N] - b[Y_s^N]\|ds + \int_0^t \|b[Y_s^N] - b[y[s]]\|ds.$$

There exists $\beta > 0$ with the following property: given $\kappa > 0$, there exists an integer $N_\kappa$, such that for each $N \geq N_\kappa$ there exists a measurable set $\Omega_N$ in the probability space on which $(Y_t^N)_{t \geq 0}$ is defined, such that $\mathbb{P}[\Omega_N] > 1 - \beta/N$, and for each sample path in $\Omega_N$ the sum of the first three terms on the right does not exceed $\kappa$ for any $t \leq u$; this statement rests upon (11), (12), and (26). Let $\lambda$ denote the Lipschitz constant of $b$, and let

$$H_t^N \equiv \sup_{0 \leq s \leq t} \|Y_s^N - y[s]\|.$$

Our inequality implies that, on $\Omega_N$,

$$H_t^N \leq \kappa + \lambda \int_0^t H_s^N ds.$$

Apply a general form of Gronwall's inequality (Ethier and Kurtz (1986), p. 498) to obtain $H_t^N \leq \kappa e^{\lambda t}$ for all $t \geq 0$. Choose $\kappa \equiv e^{-\lambda u}\delta$ to obtain the desired result (17).

### ▪ 3.2.3 Step III

[The rest of the proof takes care of details about $\sigma_N$.]. Consider the general case, where possibly $S \neq D$. On the same probability space as $(Y_t^N)_{t \geq 0}$, construct another process $(\tilde{Y}_t^N)_{t \geq 0}$ with $\tilde{Y}_0^N = Y_0^N$, which has the same rate and jump distribution at points in $S_N$, and which satisfies (12), (13) and (14) (possibly with different constants), on the whole of $I_N$. Indeed we may couple the processes, so that their trajectories actually coincide on the time interval $[0, \sigma_N]$, where $\sigma_N \equiv \inf\{t \geq 0 : Y_t^N \notin S\}$. Since $(\tilde{Y}_t^N)_{t \geq 0}$ satisfies (17) by the result of Step I, (16) follows immediately.

Next consider the case where $a \in S$. From (16) it follows that

$$\mathbb{P}[\sup_{0 \leq t \leq u} \|Y_t^N - y[t]\| \geq \delta] \leq \mathbb{P}[\sigma_N < u] + \mathbb{P}[\sup_{0 \leq t \leq u} \|\tilde{Y}_t^N - y[t]\| \geq \delta].$$

Suppose that $u < \zeta[a]$. Recall that $S$ is relatively open in $D$, and $(y[t])_{0 \leq t \leq u}$ is a continuous path which never meets the boundary of $S$. Hence there exists $\epsilon > 0$ such that the intersection with $D$ of the open $\epsilon$–tube around the path $(y[t])_{0 \leq t \leq u}$ lies within $S$. By the coupling construction, $\sigma_N$ is also $\inf\{t \geq 0 : \tilde{Y}_t^N \notin S\}$, so

$$\mathbb{P}[\sigma_N < u] \leq \mathbb{P}[\sup_{0 \leq t \leq u} \|\tilde{Y}_t^N - y[t]\| \geq \epsilon]. \tag{28}$$

Since $(\tilde{Y}_t^N)_{t \geq 0}$ satisfies (17) by the result of Step I, the last two inequalities show that $(Y_t^N)_{t \geq 0}$ also satisfies (17).

### ■ 3.2.4 Step IV

Finally consider the case where $u \equiv \zeta[a] < \infty$. By the triangle inequality,

$$\|Y_{t \wedge \sigma_N}^N - y[t]\| \leq \|Y_{t \wedge \sigma_N}^N - y[t \wedge \sigma_N]\| + \|y[t \wedge \sigma_N] - y[t]\|.$$

Using the coupling argument again,

$$\mathbb{P}[\sup_{[0,\zeta[a]]} \|Y_{t \wedge \sigma_N}^N - y[t]\| \geq 2\delta] \leq$$
$$\mathbb{P}[\sup_{[0,\zeta[a]]} \|\tilde{Y}_t^N - y[t]\| \geq \delta] + \mathbb{P}[\sup_{[0,\zeta[a]]} \|y[t \wedge \sigma_N] - y[t]\| \geq \delta].$$

We already know that the first term on the right is $O(N^{-1})$. As for the second, let

$$K \equiv \sup_{[0,\zeta[a]]} \|b[y[t]]\| < \infty.$$

In view of (27), $0 \leq r < t \leq \zeta[a]$ and $\|y[r] - y[t]\| \geq \delta$ together imply that $t - r \geq \delta/K$. Hence

$$\mathbb{P}[\sup_{[0,\zeta[a]]} \|y[t \wedge \sigma_N] - y[t]\| \geq \delta] \leq \mathbb{P}[\sigma_N \leq \zeta[a] - \delta/K].$$

Now use the reasoning of (28) to deduce that this probability too is $O(N^{-1})$. Thus (18) is established. Finally (19) follows from (18), by enlarging the state space to include time as a coördinate; the corresponding extension of the vector field $b$ is still Lipschitz. □

# 4. EXAMPLE: MULTITYPE PARTICLE SYSTEM

The following example involves only the simplest calculations, but exhibits two general techniques:

• Addition of extra components to a stochastic process to make it time–homogeneous and Markov.

• Use of an artificial time scale to simplify the solution of a differential equation.

### ■ 4.1 Quantized Particle System

Consider a model in which there are two kinds of particles — **heavy particles** and **light particles**. Moreover heavy particles may be in either of two states — **inert** or **excited**. There is also a fixed integer $w \geq 2$ which plays the following role: after $w - 1$ heavy, inert particles have been replaced by light particles, enough free energy is available to cause an inert particle to become excited.

Consider a reaction chamber containing $B$ particles, where $B$ is a random integer, divisible by $w$, with mean $\mu N$ and variance $\sigma^2 N$; the number of particles in the chamber will remain constant throughout the experiment. In the beginning, $B - 1$ particles in the chamber are heavy and inert, and one is heavy and excited.

### ■ 4.2 Dynamics of the Particle System

Here is the operation to be performed at each step:

*Select a heavy particle uniformly at random, and replace it by a light particle. Whenever the cumulative number of inert particles removed reaches a multiple of $w - 1$, some other inert particle becomes excited.*

A step at which an inert particle becomes excited is called an **excitation step**.

The number of excited particles at each step either decreases by 1 (if the particle removed was excited), or stays the same (if the particle removed was inert, not an excitation step), or increases by 1 (at an excitation step); by the same token the number of inert particles either stays the same, decreases by 2, or decreases by 1.

We may summarize the effects of the three possible types of transition by the following table:

| EVENT | $\Delta(\text{inert})$ | $\Delta(\text{excited})$ |
|---|---|---|
| particle removed is excited | 0 | $-1$ |
| particle removed is inert, excitation step | $-2$ | $+1$ |
| particle removed is inert, non–excitation step | $-1$ | 0 |

The goal is to find approximations, for large $N$, to the proportions of particles which are excited and inert, respectively, as a function of time.

### ■ 4.3 Markov Chain Notation

We model the Particle System as a process $(X_n)_{n \geq 0}$ in $\mathbb{R}^3$ as follows; we include the initial condition as one of the components so as to maintain the Markov property; we include the time index as one of the components so that the resulting process is time–homogeneous:

- $X_n^0 = B/N$ for all $n$, i.e. a rescaled number of particles.

- $X_n^1$ is the number of steps among $\{1, 2, \ldots, n\}$ at which the particle removed was inert, divided by $N$.

- $X_n^2 = n/N$, i.e. the rescaled time.

Denote by $E_n$ and $I_n$ the number of excited and inert heavy particles, respectively, after $n$ steps. We may express $E_n$ and $I_n$ in terms of $(X_n^0, X_n^1, X_n^2)$ as follows:

$$I_n = N(X_n^0 - X_n^1) - \lfloor NX_n^1/(w-1) \rfloor; \tag{29}$$

$$E_n = N(X_n^1 - X_n^2) + \lfloor NX_n^1/(w-1) \rfloor. \tag{30}$$

Given that $(X_n^0, X_n^1, X_n^2) = (x_0, x_1, x_2)$, the chance that the particle removed at step $n+1$ is inert is the proportion of heavy particles which are inert, namely

$$p_n \equiv \frac{x_0 - x_1 - N^{-1}\lfloor Nx_1/(w-1) \rfloor}{x_0 - x_2}. \tag{31}$$

We see that $(X_n)_{n \geq 0}$ is a Markov chain, where $X_n^0$ never changes, $X_n^2$ increases by $1/N$ at each step, and

$$N(X_{n+1}^1 - X_n^1) \sim \text{Bernoulli}(p_n).$$

### ■ 4.4 Mean and Covariance of Increments

In the notation of (14),

$$\mu_N[x_0, x_1, x_2] \equiv (0, N^{-1} p_n, N^{-1});$$

$$\Sigma_N[x_0, x_1, x_2] \equiv \text{diag}(0, p_n(1 - p_n)N^{-2}, 0).$$

Certainly (14) is satisfied.

### ■ 4.5 Jump Process

Let $(\nu[t])_{t \geq 0}$ be a Poisson process run such that each inter–event time $\tau_{n+1} - \tau_n$ is Exponential with mean $1/c[x]$ on the event $\{X_n = x\}$, where

$$c_N[x_0, x_1, x_2] \equiv (x_0 - x_2)N.$$

Since the rate cannot be negative, and because (31) must lie in $[0, 1]$, we choose our closed set $D \subseteq \mathbb{R}^3$ to be:

$$D \equiv \{(x_0, x_1, x_2): x_0 \geq x_1 w/(w - 1),\ x_0 \geq x_2\}.$$

As in (5), take $Y_t^N \equiv X_{\nu[t]}$, i.e. the pure jump process whose law is that of $(X_n)_{n \geq 0}$ run at the rate equal to the number of heavy particles. Thus $(Y_t)_{t \geq 0}$ stops when we run out of heavy particles. The initial condition $Y_0^N \equiv (B/N, 0, 0)$ satisfies (11) using the choice $a \equiv (\mu, 0, 0)$; this suggests that we take some constant $\kappa > \mu$, and choose the relatively open set $S \subseteq D$ to be:

$$S \equiv \{(x_0, x_1, x_2) \in D: x_0 < \kappa\}.$$

For such an $S$, it is clear that (13) holds. The product $c_N[x] \mu_N[x]$ is

$$b_N[x_0, x_1, x_2] \equiv (0, x_0 - x_1 - N^{-1} \lfloor N x_1/(w - 1) \rfloor, x_0 - x_2).$$

Hence the limiting vector field is

$$b[x_0, x_1, x_2] \equiv (0, x_0 - x_1 w/(w - 1), x_0 - x_2),$$

Since $\|b_N[x] - b[x]\| \leq 1/N$, (12) is trivially satisfied. Moreover $b$ is linear, hence globally Lipschitz on $D$. All the assumptions of the Theorem have now been validated.

### ■ 4.6 Fluid Limit

The ODE (15) takes the form:

$$(\dot{y}_0[t], \dot{y}_1[t], \dot{y}_2[t]) = (0, y_0 - y_1 w/(w - 1), y_0 - y_2);\ y[0] = (\mu, 0, 0).$$

Thus $y_0[t] = \mu$ for all $t$, while the other components satisfy separate linear ODE's:

$$\frac{dy_1}{\mu - y_1 w/(w-1)} = dt; \quad \frac{dy_2}{\mu - y_2} = dt,$$

with solution

$$y_1[t] = \frac{w-1}{w}\mu(1 - e^{-tw/(w-1)}); \quad y_2[t] = \mu(1 - e^{-t}).$$

Hence the fluid limit of the number of heavy particles, divided by $N$, is

$$h[t] \equiv (y_0[t] - y_2[t]) = \mu e^{-t}.$$

It follows from (29) and (30) that the fluid limits of the numbers of excited and inert particles, respectively, divided by $N$, are $e[t]$ and $h[t] - e[t]$, where

$$e[t] \equiv y_1[t]w/(w-1) - y_2[t] = \mu(e^{-t} - e^{-tw/(w-1)}).$$

$$= h[t] - \mu(h[t]/\mu)^{w/(w-1)}.$$

### ■ 4.7 Conclusion

Since one heavy particle is replaced by a light particle at each step, it is natural to express the fluid limit in terms of $h$, the number of heavy particles divided by $N$, which decreases linearly from $\mu$ (at time 0) to 0 (at time $\mu$), where "time" means number of steps, divided by $N$. In those terms, $h - \mu(h/\mu)^{w/(w-1)}$ and $\mu(h/\mu)^{w/(w-1)}$ are the fluid limits of the numbers of excited and inert particles, divided by $N$, respectively; uniform convergence in probability to these limits occurs in the sense of the Theorem.

**Acknowledgments:** J. R. Norris (Cambridge University) helped to clarify various ideas. Thanks to Warwick de Launay and David P. Moulton for expressing interest in seeing a presentation of this material.